\documentclass{article}
\usepackage{amsfonts}
\title{\large \bf On equality of central and class preserving automorphisms of finite $p$-groups }
\author{\small \bf Hemant Kalra and Deepak Gumber \\
\small \em School of Mathematics and Computer Applications\\
\small \em Thapar University, Patiala - 147 004,
India\\
\small E-mail: dkgumber@yahoo.com, happykalra26@gmail.com}

\date{}
\newtheorem{thm}{Theorem}[section]
\newtheorem{lm}[thm]{Lemma}

\begin{document}

\maketitle
\begin{abstract}
Let $G$ be a finite non-abelian $p$-group, where $p$ is a prime. Let $\mathrm{Aut}_c(G)$ and
$\mathrm{Aut}_z(G)$ respectively denote the group of all class preserving and
central automorphisms of $G$. We give a necessary and sufficient condition for $G$ such that $\mathrm{Aut}_c(G)=\mathrm{Aut}_z(G)$ and classify all finite non-abelian $p$-groups $G$ with elementary abelian or cyclic center  such that $\mathrm{Aut}_c(G)=\mathrm{Aut}_z(G).$ We also  characterize all finite $p$-groups $G$ of order $\leq p^7$ such that $\mathrm{Aut}_c(G)=\mathrm{Aut}_z(G)$ and complete the classification of all finite $p$-groups of order $\le p^5$ for which there exist non-inner class preserving automorphisms.
\end{abstract}

\vspace{2ex}

\noindent {\bf 2010 Mathematics Subject Classification:} 20D15,
20D45.

\vspace{2ex}

\noindent {\bf Keywords:} Class preserving automorphism, Central
automorphism, Camina pair.

\section{\large Introduction}
Let $G$ be a finite group. Let $\mathrm{Aut}(G)$  and
$\mathrm{Inn}(G)$ respectively denote the group of all automorphisms
and inner automorphisms of $G$. An automorphism $\alpha$ of $G$ is
called a class preserving automorphism if $\alpha(x)\in x^G$, the
conjugacy class of $x$ in $G$, for all $x\in G$; and an automorphism
$\phi$ of  $G$ is called a central automorphism if $g^{-1}\phi(g)\in
Z(G)$ for all $g\in G$. The set $\mathrm{Aut}_c(G)$ of all class
preserving automorphisms and the set $\mathrm{Aut}_z(G)$ of all
central automorphisms of $G$ are normal subgroups of $\mathrm{Aut}
(G)$. Let $\mathrm{Aut}_{z}^{z}(G)$ denote the group of all central
automorphisms of $G$ fixing the center $Z(G)$ of $G$ element-wise.
Interest in the equality of class preserving and inner automorphisms
dates back to 1911 when Burnside [4, p. 463] posed the following
question: Does there exist a finite group $G$ such that $G$ has a
non-inner class preserving automorphism?  For details, readers can
see the excellent survey article by Yadav [19]. In the recent past,
interest of many mathematicians turned on the equalities of
$\mathrm{Aut}_z(G)$ and $\mathrm{Inn}(G)$, $\mathrm{Aut}_z(G)$ and
$Z(\mathrm{Inn}(G))$, $\mathrm{Aut}_z(G)$ and
$\mathrm{Aut}_{z}^{z}(G)$, and $\mathrm{Aut}_{z}^{z}(G)$ and
$\mathrm{Inn}(G)$ (see for example [2], [5], [6], [7], [10] and
[18]). Curran and McCaughan [6] characterized finite $p$-groups $G$
for which $\mathrm{Aut}_z(G)=\mathrm{Inn}(G)$. They proved that if
$G$ is a finite $p$-group, then $\mathrm{Aut}_z(G)=\mathrm{Inn}(G)$
if and only if $\gamma_{2}(G)=Z(G)$ and $Z(G)$ is cyclic. As a
consequence of our results, we give an easy and very short proof of
this result in section 3. Yadav [19], in his survey article, has
asked to classify all finite $p$-groups $G$ of class $2$ such that
$\mathrm{Aut}_c(G)=\mathrm{Aut}_z(G)$. In section 3, we give a
necessary and sufficient condition for a finite non-abelian $p$-group $G$ such that
$\mathrm{Aut}_c(G)=\mathrm{Aut}_z(G)$. We also classify all finite $p$ groups $G$ of class 2 for which $\mathrm{Aut}_c(G)=\mathrm{Aut}_z(G)$ in the cases when $Z(G)$
is cyclic or elementary abelian. In section 4, we characterize all finite $p$-groups
$G$ of order $\leq p^7$ for which
$\mathrm{Aut}_c(G)=\mathrm{Aut}_z(G)$. We also complete the
classification of all finite $p$-groups $G$ of order $\le p^5$ for
which $\mathrm{Aut}_c(G)\ne\mathrm{Inn}(G)$. It follows from [12]
that for extra special $p$-groups, and hence for all groups of order
$p^3$, all class preserving automorphisms are inner. From [13], it
follows that  $\mathrm{Aut}_c(G)=\mathrm{Inn}(G)$ for all groups $G$
of order $p^4$. Yadav [17] has found all the groups of order $p^5$,
$p$ an odd prime, for which $\mathrm{Aut}_c(G)\ne\mathrm{Inn}(G)$. In
the present paper we show that out of $51$ groups of order 32, there
are only two groups for which $\mathrm{Aut}_c(G)\ne\mathrm{Inn}(G)$.
In section 5, as an application of our results, we find all finite
$p$-groups $G$ of order $p^n(n\leq 6)$ for which
$\mathrm{Aut}_c(G)=\mathrm{Aut}_z(G)$ from the lists of all such
$p$-groups given by James [11] (for odd $p$) and Hall and
Senior [8] (for $p=2$).

\section{\large Notations and Preliminaries}
By $\mathrm{Hom}(G,A)$ we denote the group of all homomorphisms of
$G$ into an abelian group $A$, by $d(G)$ we denote the rank of $G$,
and by $C_{p^n}$ we denote the cyclic group of order $p^n$. A
non-abelian group $G$ that has no non-trivial abelian direct factor
is said to be purely non-abelian. Observe that a group $G$ is purely
non-abelian if its center Z(G) is contained in the frattini subgroup
$\Phi(G)$ of $G$. For two subgroups $H$ and $K$ of $G$, $[H,K]$
denotes the subgroup of $G$ generated by all commutators
$[x,y]=x^{-1}y^{-1}xy$ with $x\in H$ and $y\in K$. By $[x,G]$ we
denote the set of all commutators of the form $[x,g],\;g\in G$.
Observe that if $G$ is nilpotent of class 2, then $[x,G]$ is a
normal subgroup of $G$. The lower central series of a group $G$ is
the descending series
$G=\gamma_1(G)\geq\gamma_2(G)\geq\cdots\geq\gamma_i(G)\geq\cdots$,
where $\gamma_{n+1}(G)=[\gamma_n(G),G]$, and upper central series is
the ascending series $Z(G)=Z_1(G)\leq Z_2(G)\leq\cdots \leq
Z_i(G)\leq\cdots$, where $Z_{i+1}(G)=\{x\in G|[x,y]\in Z_i(G)\;
\mathrm{for\; all}\; y\in G\}.$ Let $G$ be a finite $p$-group and
$M$ be non-trivial proper normal subgroup of $G$. The pair $(G,M)$
is called a Camina pair if $M\subseteq [y,G]$ for all $y\in G-M$ and
$G$ is called a Camina group if
$(G,\gamma_{2}(G))$ is a Camina pair.\\

The following well known result of Adney and Yen [1] will be
referred to as Adney-Yen Lemma.
\begin{lm}
If $G$ is a purely non-abelian group, then there is a ono-to-one correspondence between $\mathrm{Aut}_z(G)$ and $\mathrm{Hom}(G/\gamma_{2}(G),Z(G))$.
\end{lm}

And the following result of Morigi [15, Lemma 0.4] will be referred
to as Morigi Lemma.

\begin{lm}
Let $G$ be a finite nilpotent group of class $2$. Then $\exp(\gamma_2(G))=\exp(G/Z(G))$ and in the decomposition of $G/Z(G)$ in direct product of cyclic groups at least two factors of maximal order must occur.
\end{lm}

The well known commutator identities
$$[x, yz] = [x, z][x, y][x, y, z]\;;\;\;[xy, z] = [x, z][x, z, y][y, z],$$
where $x, y, z \in G$, will be frequently used without any reference.

\section{\large Main Results}
We start with the following theorem which gives a necessary and sufficient condition on a finite $p$-group $G$ such that $\mathrm{Aut}_c(G)=\mathrm{Aut}_{z}(G)$.

\begin{thm} Let $G$ be a finite $p$-group. Then
$\mathrm{Aut}_c(G)=\mathrm{Aut}_{z}(G)$ if and only if
$\mathrm{Aut}_c(G)\approx \mathrm{Hom}(G/Z(G),\gamma_{2}(G))$ and
$\gamma_2(G)=Z(G)$.
\end{thm}

\noindent {\bf Proof.} First assume that $\mathrm{Aut}_c(G)=\mathrm{Aut}_{z}(G)$. We show that $Z(G)\leq \Phi(G).$  Assume contrarily that $Z(G)$ is not contained in $\Phi(G)$. Choose an element $h\in Z(G)-M$ for some maximal subgroup $M$ of $G$. Therefore, $G=M\langle h\rangle$. There exists a non-trivial element $z$ in $Z(G)\cap \Phi(G)$ of order $p$. It is easy to check that the map $\mu :G\rightarrow G$ defined by $\mu(mh^i)=mh^iz^i$ for every $m\in M$ and for every $i,\;0\leq i\leq p-1,$ is a central automorphism of $G$ which is not class preserving and thus $Z(G)\leq \Phi(G).$  For any commutator $[a,b]\in G$, where $a,b\in G$, we can define an inner automorphism $\iota_b:G\rightarrow G$ such that $a^{-1}\iota_b(a)=[a,b]\in Z(G)$. Thus $\gamma_2(G)\le Z(G)$.
For any $\mu \in
\mathrm{Aut}_{z}(G)$, the map
$\psi_{\mu}:G/Z(G)\rightarrow \gamma_{2}(G)$ defined as 
$\psi_{\mu}(bZ(G))=b^{-1}\mu (b)$ is a homomorphism. It is easy to see that the map $\psi$ sending $\mu$ to $\psi_{\mu}$ is a monomorphism of the group $\mathrm{Aut}_z(G)$ into the group $\mathrm{Hom}(G/Z(G),\gamma_{2}(G))$. For any $\tau \in
\mathrm{Hom}(G/Z(G),\gamma_{2}(G))$, the map $\mu:G\rightarrow G$
defined as $\mu(g)=g\tau (gZ(G))$, where $g\in G$, is a central automorphism and
$\psi(\mu)=\psi_{\mu}=\tau$. Thus $\psi$ is onto as well and hence $\mathrm{Aut}_c(G)\approx \mathrm{Hom}(G/Z(G),\gamma_{2}(G))$.
 But, by Adney-Yen Lemma, $|\mathrm{Aut}_z(G)|=|\mathrm{Aut}_c(G)|=|\mathrm{Hom}(G/\gamma_{2}(G),Z(G))|$, and so
\begin{equation}
|\mathrm{Hom}(G/Z(G),\gamma_{2}(G))|=
|\mathrm{Hom}(G/\gamma_{2}(G),Z(G))|.
\end{equation}
Suppose to the contrary that $\gamma_2(G)<Z(G)$. Then $G/Z(G)$ is a
proper quotient of $G/\gamma_2(G)$ and $|(G/\gamma_2(G))/(G/Z(G))|
=|Z(G)/\gamma_2(G)|>1.$ It thus follows from [5, Lemma 2.8] that
$\mathrm{Hom}(G/Z(G),\gamma_{2}(G))$ is isomorphic to a proper
subgroup of $\mathrm{Hom}(G/\gamma_2(G),Z(G))$. This is a
contradiction to (1) and hence $\gamma_2(G)=Z(G)$. 

Conversely assume that $\mathrm{Aut}_c(G)\approx \mathrm{Hom}(G/Z(G),\gamma_{2}(G))$ and
$\gamma_2(G)=Z(G)$. Observe that since $\gamma_2(G)=Z(G)$, $\mathrm{Aut}_c(G)\le\mathrm{Aut}_z(G)$ and $G$ is purely non-abelian. By Adney-Yen Lemma 
$$|\mathrm{Aut}_z(G)|=|\mathrm{Hom}(G/\gamma_2(G),Z(G))|=|\mathrm{Hom}(G/Z(G),\gamma_2(G))|=|\mathrm{Aut}_c(G)|.$$
This completes the proof.
\hfill $\Box$\\

As a consequence of Theorem 3.1, we give the following easy proof of
main result of Curran and McCaughan [6]. The proof of only if part
is easy and is similar to as given in [6]. But we give it here for the sake of completeness.

\begin{thm}
If $G$ is a finite $p$-group,
then $\mathrm{Aut}_z(G)=\mathrm{Inn}(G)$ if and only if
$\gamma_{2}(G)=Z(G)$ and $Z(G)$ is cyclic.
\end{thm}
\noindent {\bf Proof.} Suppose $\gamma_{2}(G)=Z(G)$ and $Z(G)$ is cyclic. Then $\mathrm{Inn}(G)\leq \mathrm{Aut}_z(G),$ $\exp(G/Z(G))=\exp(\gamma_2(G))$ by Morigi Lemma, and $G$ is purely non-abelian. By Adney-Yen lemma
$$|\mathrm{Aut}_z(G)|=
|\mathrm{Hom}(G/\gamma_2(G),Z(G))|=|G/Z(G)|=
|\mathrm{Inn}(G)|,$$ and thus $\mathrm{Aut}_z(G)=\mathrm{Inn}(G)$.

Conversely suppose that $\mathrm{Inn}(G)=\mathrm{Aut}_z(G)$. Then, as in above theorem, we can prove that   $\mathrm{Aut}_z(G)\approx
\mathrm{Hom}(G/Z(G),\gamma_2(G))$ and $\gamma_2(G)=Z(G)$. Thus nilpotency class of $G$ is 2 and hence $\exp(G/Z(G))=\exp(\gamma_2(G))$. It now follows from $$\mathrm{Hom}(G/Z(G),\gamma_2(G))\approx \mathrm{Aut}_z(G)=\mathrm{Inn}(G)\approx G/Z(G),$$ that $\gamma_2(G)$ is cyclic.
\hfill $\Box$

\begin{thm}
Let $G$ be a finite non-abelian $p$-group such that $Z(G)$ is elementary abelian. Then $\mathrm{Aut}_c(G)=\mathrm{Aut}_z(G)$ if and only if $G$ is a Camina $p$-group of nilpotency class $2$.
\end{thm}
\noindent {\bf Proof.} Suppose first that $G$ is a Camina $p$-group
of nilpotency class 2. Then $\gamma_2(G)=Z(G)$ by [14, Lemma 2.1]. Let $\phi \in \mathrm{Aut}_c(G)$ and $g\in
G$. Then $g^{-1}\phi(g)\in \gamma_2(G)= Z(G)$. Thus
$\phi$ is a central automorphism. On the other hand, let $\psi\in
\mathrm{Aut}_z(G)$ and let $h\in G$.  If $h\in Z(G)$, then $\psi(h)=h$ and if $h\in G-Z(G)$, then $h^{-1}\psi (h)\in
Z(G)\subseteq [h,G]$. It thus follows that $\psi(h)=b^{-1}hb$ for some
$b\in G$.  Therefore $\psi$ is a class preserving automorphism and hence
$\mathrm{Aut}_c(G)=\mathrm{Aut}_z(G)$.

Conversely suppose that $\mathrm{Aut}_c(G)=\mathrm{Aut}_z(G)$. By Theorem 3.1,
$\gamma_2(G)=Z(G)$ and
$\mathrm{Aut}_c(G)\approx\mathrm{Hom}(G/Z(G),\gamma_{2}(G)).$   It follows that
$\mathrm{exp}(G/Z(G))=
\mathrm{exp}(\gamma_{2}(G))=\mathrm{exp}(Z(G))=p$.
Let $|G|=p^{l}$,
$\{y_{1},y_{2}, \ldots,y_{r}\}$ be a minimal generating set of
$G$ and let $|\gamma_{2}(G)|=p^{s}$. Then
$$|\mathrm{Aut}_c(G)|=|\mathrm{Hom}(G/Z(G),\gamma_{2}(G))|
=p^{s(l-s)}=p^{sr}.$$
Any element $y\in G-\gamma_{2}(G)$ is a part of a
minimal generating set $\{y=y_{1},y_{2},\ldots ,y_{r}\}$ for $G$.
If possible suppose $[y,G]< \gamma_{2}(G)$, then $|y^{G}|=|[y,G]|<
|\gamma_{2}(G)|=p^s$. Thus $|\mathrm{Aut}_c(G)|< p^{sr}$,
a contradiction and hence $G$ is a Camina $p$-group of class 2.\hfill $\Box$

\begin{thm}
Let $G$ be a finite non-abelian $p$-group such that $Z(G)$ is cyclic. Then $\mathrm{Aut}_c(G)=\mathrm{Aut}_z(G)$ if and only if $Z(G)=\gamma_2(G).$
\end{thm}
\noindent {\bf Proof.} If $Z(G)=\gamma_2(G)$, then  $\mathrm{Inn}(G)=\mathrm{Aut}_z(G)$ by Theorem 3.2. Let $\mu$ be any class preserving automorphism of $G$ and let $g\in G$. Then $g^{-1}\mu(g)\in \gamma_2(G)=Z(G)$. Thus $\mu$ is a central automorphism and hence $$\mathrm{Aut}_c(G)\leq \mathrm{Aut}_z(G)=\mathrm{Inn}(G)\leq \mathrm{Aut}_c(G).$$
Converse follows from Theorem 3.1.\hfill $\Box$

\section{\large Groups of order $p^n,\;n\leq 7$}
In this section, we  characterize all finite $p$-groups of order
$p^n(n\leq 7)$ such that $\mathrm{Aut}_c(G)=\mathrm{Aut}_z(G)$. For
this, we need the following theorem which also completes the
classification of  all finite $p$-groups of order $\le p^5$ such that $\mathrm{Aut}_c(G)\ne\mathrm{Inn}(G)$. In the
theorem, we show that out of 51 groups of order 32, there are only 2
groups for which $\mathrm{Aut}_c(G)\ne \mathrm{Inn}(G)$. A list of
groups of order 32 is available from  Hall and Senior [8] and from
Sag and Wamsley [16] with minimal presentations of groups. In [8],
groups of order 32 are divided into 8 isoclinism families. We denote
the $i$th family by $\Phi_i$. As mentioned in [16], the groups are
listed in same order in both [8] and [16] and so we take the liberty
of choosing the presentation of a group from any of the lists. Unless or otherwise stated, we use
Sag and Wamsley's list and adopt the same notations for the
nomenclature and presentations of the groups. However, we write the
$i$-th group of order 32 in the list as $G_i$ and generators 1,2,3
and 4  respectively as $x,y,z$ and $w$. Yadav [17] has proved that
if $G$ and $H$ are two finite non-abelian isoclinic groups, then
$\mathrm{Aut}_c(G)\approx \mathrm{Aut}_c(H)$. Therefore, it is sufficient
to pick only one group from each isoclinism family. We need the
following lemma to prove the theorem.

\begin{lm}
If $G$ is a group of order $p^n$, $n\geq 3$ and $|Z(G)|=p^{n-2}$,
then $\mathrm{Aut}_c(G)=\mathrm{Inn}(G)$.
\end{lm}
{\bf Proof.} Observe that nilpotency class of $G$ is $2$. If $x$ is
a non-central element of $G$, then $|C_G(x)|=p^{n-1}$ and thus
$|x^G|=|[x,G]|=p$. If $x$ is a central element of $G$, then
$|x^G|=|[x,G]|=1$. In any case, $[x,G]$ is cyclic and therefore
$\mathrm{Aut}_c(G)=\mathrm{Inn}(G)$ by [17, Theorem 3.5]. \hfill
$\Box$

\begin{thm}
 For all groups G of order $32$, except for the groups
$G_{44}$ and
$G_{45}$ in the sixth family $\Phi_6$, $\mathrm{Aut}_c(G)=\mathrm{Inn}(G)$.
\end{thm}
{\bf Proof.} The first family $\Phi_1$ contains abelian groups $G_{1}$ to $G_{7}$ and therefore all
class preserving automorphisms for these groups are inner. From the second family $\Phi_2$, we pick $G_{11}=\{x^{2}, y^{4}, z^{2}, [x,y,x], [x,y,y], [x,z], [y,z]\}$. Clearly $z$ is in the center, and since $[x,y^2]=[x,y]^2=[x^2,y]=1$,
it follows that $y^2$ is in the center and $|[x,y]|=2$. Also, since
$z$ commutes with $x$ and $y$, it commutes with $[x,y]$ as well.
Thus $[x,y]$ is also in the center and hence $|Z(G_{11})|=8$. Therefore
$\mathrm{Aut}_c(G_{11})=\mathrm{Inn}(G_{11})$ by above lemma.

From third isoclinism family, we take $G_{23}=\{x^{8}, y^{2}, z^{2}, [x,y]x^2,
[x,z], [y,z]\}$. Observe that $G_{23}=H\oplus K,$ where $H=\{x,y| x^8=y^2=1, yx=x^7y\}$
and $K=\{z|z^2=1\}$. It is easily seen that
 $K$ is cyclic  and therefore $\mathrm{Aut}_c(K)=\mathrm{Inn}(K)$. Also,
since every element of $H$ can be written as $x^iy^j$ for suitable $i$ and
$j$, $\mathrm{Aut}_c(H)=\mathrm{Inn}(H)$ by [13,
Proposition 4.1]. Hence
$\mathrm{Aut}_c(G_{23})=\mathrm{Inn}(G_{23})$ by [12, Proposition
2.2].

In $\Phi_4$, consider  $G_{34}=\{x^{4}, y^{4}, z^{2}, [x,y],
[x,z]x^{2}, [y,z]y^{2}\}$. Let $H_{34}=\langle x,y\rangle$. Then $H_{34}$ is an
abelian normal subgroup of order 16, and therefore $G_{34}/H_{34}$ is
a cyclic group of order 2. Hence
$\mathrm{Aut}_c(G_{34})=\mathrm{Inn}(G_{34})$ by [9, Proposition
2.7].

From fifth family, consider $$G_{43}=\{w^2, x^2y^{-2}, x^2z^{-2},
[z,y]x^2, [w,x]x^2, [x,y], [x,z], [y,w], [z,w]\}.$$ Since
$d(G_{43})=4$, $|\Phi(G_{43})|=|\gamma_2(G_{43})|=2$. Also, since
$G_{43}$ is a stem group, $Z(G_{43})\le \gamma_2(G_{43})$. Thus
$G_{43}$ is an extra-special group and hence
$\mathrm{Aut}_c(G_{43})=\mathrm{Inn}(G_{43})$ by [12, Theorem 3.2].

Next consider the group $G=G_{44}=\{x^{2}, z^{2}, [y,x]y^{4},
[y,z]y^{2}, [x,z]\}$ from $\Phi_6$.  Since $yx=xy[y,x]=xyy^{-4}=xy^{-3}$, and since $[y,z]=y^{-2}$ implies that $zy=y^{-1}z$, it
follows that every element of $G$ can be written in the
form $x^iy^jz^k$, where $0\leq i,k\leq 1$, and $0\leq j\leq 7$; and hence $|y|=8$. Since $G$ is of nilpotency class
3 and $d(G)=3$, it follows that $|\Phi(G)|=|\gamma_2(G)|=4$. But $y^2=[z,y]\in \gamma_2(G)$ is of order 4, therefore $\gamma_2(G)=\langle y^2\rangle$. Now
$G_{44}$ is a stem group, therefore $Z(G)< \gamma_2(G)$. Since $[z,y^4]=[z,y^2]^2=[z,y]^4=y^8=1$
and $[x,y^2]=[x,y]^2=y^8=1$, $y^4\in Z(G)$ and thus
$Z(G)=\langle y^4\rangle.$ Therefore $G/Z(G)$ is a class 2 group of order 16.
Then $Z(G/Z(G))$ is of order 4 and hence $|Z_2(G)|=8$. We prove that $Z(G)\subseteq [g,G]$ for
all $g\in G-\gamma_2(G)$.
Let $g=x^iy^jz^k\in G-\gamma_2(G)$. First suppose that $j$ is even. Then $g=y^jx^iz^k$. Both $i$ and $k$ cannot be zero, because then $g\in \gamma_2(G)$. If $k=1$, then $[g,y^2]=[y^jx^iz,y^2]=[x^iz,y^2]=[z,y^2]=y^4$. If $k=0$, then $[g,y]=[y^jx,y]=[x,y]=y^4$. Thus $Z(G)\subseteq [g,G]$ for
all $g\in G-\gamma_2(G)$ in this case. Next suppose that $j$ is odd. Then $g=y^{j-1}x^iyz^k$. If $k=1$, then $[g,y^2]=[y^{j-1}(x^iyz),y^2]=[x^i(yz),y^2]=[yz,y^2]=[z,y^2]
=[z,y]^2=y^4$. Suppose $k=0$ and $i=1$. Then $[g,y]=[y^{j-1}(xy),y]=[xy,y]=[x,y]=y^4$. Finally suppose $k=0=i$. Then $[g,x]=[y^{j-1}y,x]=[y,x]=y^4$. Thus  $Z(G)\subseteq [g,G]$ for
all $g\in G-\gamma_2(G)$ in this case as well.
 Therefore, by [17, Lemma 2.2]
$$|\mbox{Aut}_c(G)|\geq |\mbox{Aut}_z(G)||G/Z_2(G)|.$$
Since $Z(G)<\gamma_2(G)$, $G$ is purely non abelian and therefore by
Adney-Yen Lemma, we have
$$|\mbox{Aut}_z(G)|=|\mbox{Hom}(G/\gamma_2(G),Z(G))|=|\mbox{Hom}(C_2\times C_2\times C_2,C_2)|=8.$$
 Thus  $|\mbox{Aut}_c(G)|\geq
2^5>2^4=|\mbox{Inn}(G)|$.

From seventh family, we take $G_{46}=\{x^{2}, y^{4}, [x,y,x],
[x,y,y,y]\}$. Let $u=[x,y]$ and $v=[x,y,y]=[u,y]$. Then $u$ commutes
with $x$ and $v$ commutes with $y$. Observe that
$1=[x^2,y]=[x,y]^2=u^2\;;\;[x,y^2]=[x,y]^2[x,y,y]=v\;
;\;1=[x,y^4]=[x,y^2]^2[x,y^2,y^2]=v^2.$ Thus $|u|=|v|=2$. Since
$xy^2 = (xy)y = (yxu)y = yx(uy) = yx(yuv) = y(xy)uv = y(yxu)uv =
y^2xv$, $vx = (xy^2xy^2)x = (xy^2x)y^2x = (xy^2x)(xy^2v) = xv$. Thus
$v$ commutes with $x$ and hence with $u$ as well. It is easy to see
that every element of $G_{46}$ can be written in the form
$x^iy^ju^kv^l$, where $0\leq i,k,l\leq 1$ and $0\leq j\leq 3$. Take
$H_{46}=\langle x,u,v\rangle$. Then $H_{46}$ is an abelian normal subgroup of
order 8 such that $G_{46}/H_{46}=\langle yH_{46}\rangle$ is a cyclic group of
order 4. Thus $\mathrm{Aut}_c(G_{46})=\mathrm{Inn}(G_{46})$ by [9,
Proposition 2.7].

Finally consider $G_{51}=\{x^{8}y^{-2}, [y,x^{7}]x^{2}\}$ from
$\Phi_8$. The relation $[x^7,y]=x^2$ implies that $yx^9=x^7y$. Post
multiplying by $x^8=y^2$ we have $yx=x^7y^3=x^{15}y$. Thus every element of
$G_{51}$ can be written as $x^iy^j,\;0\le i\le 15,0\le j\le 1$ and
hence $\mathrm{Aut}_c(G_{51})=\mathrm{Inn}(G_{51})$ by [13,
Proposition 4.1].
\hfill $\Box$ \\

\noindent {\bf Remark:} The order of $\mbox{Aut}_c(G)$ for
$G=G_{44}$ in the above theorem is in fact exactly equal to 32.
Observe that $C_G(x)=\langle x,y^2,z\rangle,\;C_G(y)=\langle y\rangle$ and
$C_G(z)=\langle x,y^4,z\rangle$. Thus $|C_G(x)|=16$ and $|C_G(y)|=|C_G(z)|=8$.
Hence $|x^G|=2$ and $|y^G|=|z^G|=4$. Since any class preserving
automorphism preserves the conjugacy classes, there are
$|x^G|,\;|y^G|$ and $|z^G|$ choices for the images of  $x,\;y$ and
$z$ respectively under it. Thus $|\mbox{Aut}_c(G)|\leq
|x^G||y^G||z^G|=32$.

\begin{thm}
Let $G$ be a non-abelian group of order $p^n,\;3\leq n\leq 5.$ Then $\mathrm{Aut}_c(G)=\mathrm{Aut}_z(G)$ if and only if $\gamma_2(G)=Z(G)$ and $Z(G)$ is cyclic.
\end{thm}
\noindent {\bf Proof.} Suppose
$\mathrm{Aut}_c(G)=\mathrm{Aut}_z(G)$. Then $\gamma_2(G)=Z(G)$ by
Theorem 3.1. If $|G|=p^3$ or $p^4$, then
$\mathrm{Aut}_c(G)=\mathrm{Inn}(G)$ by [12, 13]. If $|G|=p^5$, then
since class of $G$ is 2, it follows from Theorem 4.2 (for $p=2$) and from
[17, Theorem 5.5] (for odd $p$) that
$\mathrm{Aut}_c(G)=\mathrm{Inn}(G).$ Thus
$\mathrm{Inn}(G)=\mathrm{Aut}_z(G)$ and hence $\gamma_{2}(G)$ is cyclic by
Theorem 3.2. The converse follows from Theorem 3.4. \hfill $\Box$

\begin{thm}
Let $G$ be a non-abelian group of order $p^6$. Then $\mathrm{Aut}_c(G)=\mathrm{Aut}_z(G)$ if and only if either $\gamma_2(G)=Z(G)$ and $Z(G)$ is cyclic or $G$ is a Camina $p$-group of nilpotency class $2$.
\end{thm}
\noindent {\bf Proof.} Suppose $\mathrm{Aut}_c(G)=\mathrm{Aut}_z(G)$. Then  $\gamma_2(G)=Z(G)$ by Theorem 3.1. If $Z(G)$ is not cyclic, then we prove that $Z(G)$ is elementary abelian and the result will then follow from Theorem 3.3. Now $p^2\le |Z(G)|\le p^4$. If $|Z(G)|=p^2$, then it is trivially elementary abelian. Let $|Z(G)|=p^3$ and  let $Z(G)\approx C_{p^2}\times C_p$. Then $\exp(Z(G))=\exp(\gamma_2(G))=\exp(G/Z(G))=p^2$ and hence $G/Z(G)\approx C_{p^2}\times C_p$, a contradiction to Morigi Lemma. If $|Z(G)|=p^4$, then $G/Z(G)$ is elementary abelian. Therefore $\exp(Z(G))=\exp(\gamma_2(G))=\exp(G/Z(G))=p$ and hence $Z(G)$ is elementary abelian. The converse follows from Theorems 3.3 and 3.4.
\hfill $\Box$

\begin{thm}
Let $G$ be a non-abelian group of order $p^7$. Then $\mathrm{Aut}_c(G)=\mathrm{Aut}_z(G)$ if and only if either $\gamma_2(G)=Z(G)$ and $Z(G)$ is cyclic or $G$ is a Camina $p$-group of nilpotency class $2$.
\end{thm}
\noindent {\bf Proof.} Suppose $\mathrm{Aut}_c(G)=\mathrm{Aut}_z(G)$. Then $\gamma_2(G)=Z(G)$ by Theorem 3.1.
If $Z(G)$ is not cyclic, then as in Theorem
4.4, we need only prove that $Z(G)$ is elementary abelian. Now
$p^2\le |Z(G)|\le p^5$. The cases  $|Z(G)|=p^2$ or $p^5$ can be
handled as in the above theorem. Let $|Z(G)|=p^3$ and  let
$Z(G)\approx C_{p^2}\times C_p$. Then $G/Z(G)\approx C_{p^2}\times
C_{p^2}$ by Morigi Lemma. Thus $G$ is a 2-generator class 2 group
and hence $\gamma_2(G)$ is cyclic by [3, Lemma 36.5]. This is a
contradiction to $\gamma_2(G)=Z(G)\approx C_{p^2}\times C_p$. Let
$|Z(G)|=p^4$ and let $\exp(Z(G))=p^2$ or $p^3$. Then $|G/Z(G)|=p^3$
and $\exp(G/Z(G))=p^2$ or $p^3$, which is not possible by Morigi
Lemma.  The converse follows from Theorems 3.3 and 3.4.
\hfill $\Box$

\section{\large Application}

In this section, we use the classification of all groups of order
$p^n$, $ 5\le n\le 6$, given by James [11] for odd $p$ and by Hall and
Senior [8] for $p=2$. As an application of our results, we find
those groups $G$ of order $p^5$ and $p^6$ for which
$\mathrm{Aut}_c(G)=\mathrm{Aut}_z(G)$.

\begin{thm}
Let $G$ be a non-abelian group of order $p^5$, where $p$ is an odd prime. Then $\mathrm{Aut}_c(G)=\mathrm{Aut}_z(G)$ if and only if $G\in \Phi_{5}$.
\end{thm}
\noindent {\bf Proof.} By Theorem 4.3,
$\mathrm{Aut}_c(G)=\mathrm{Aut}_z(G)$ if and only if
$\gamma_2(G)=Z(G)$ and $Z(G)$ is cyclic. This happens if and only if
$G\in \Phi_5$ by $\S 4.1$ of [11].
\hfill $\Box$

\begin{lm}
A non-abelian group $G$ of order $p^6$, where $p$ is an odd prime, is a Camina group of class $2$ if and only if $G\in \Phi_{15}$.
\end{lm}
\noindent {\bf Proof.} In a Camina $p$-group of class 2, $\gamma_2(G)=Z(G)$ and $|x^G|=|\gamma_2(G)|$ for each $x\in G-Z(G)$. This happens if and only if $G\in \Phi_{15}$ by $\S 4.1$ of [11].
\hfill $\Box$

\begin{thm}
Let $G$ be a non-abelian group of order $p^6$, where $p$ is an odd prime. Then $\mathrm{Aut}_c(G)=\mathrm{Aut}_z(G)$ if and only if $G\in \Phi_{14}$ or $\Phi_{15}$.
\end{thm}
\noindent {\bf Proof.} By Theorem 4.4,
$\mathrm{Aut}_c(G)=\mathrm{Aut}_z(G)$ if and only if either $G$ is a
Camina $p$-group of nilpotency class 2 or  $\gamma_2(G)=Z(G)$ and
$Z(G)$ is cyclic. Now  $G$ is a
Camina $p$-group of class 2 if and only if $G\in \Phi_{15}$ by above lemma, and  $\gamma_2(G)=Z(G)$ and $Z(G)$ is cyclic if
and only if $G\in \Phi_{14}$ by $\S 4.1$ of [11].
\hfill $\Box$

\begin{thm}
Let $G$ be a non-abelian group of order $2^5$. Then $\mathrm{Aut}_c(G)=\mathrm{Aut}_z(G)$ if and only if  $G\in \Phi_{5}$.
\end{thm}
\noindent {\bf Proof.} By Theorem 4.3, $\mathrm{Aut}_c(G)=\mathrm{Aut}_z(G)$ if and only if $\gamma_2(G)=Z(G)$ and $Z(G)$ is cyclic. There are only two families $viz. \;\Phi_4$ and $\Phi_5$ which consist of groups $G$ such that $\gamma_2(G)=Z(G)$. Consider the group $H=G_{34}=\{x^{4}, y^{4}, z^{2}, [x,y], [x,z]x^{2},
[y,z]y^{2}\}$ from $\Phi_4$. Observe that $[x^2,z]=[x,z]^2=x^4=1$ and $[y^2,z]=[y,z]^2=y^4=1$. Thus $x^2,y^2\in Z(H)=\gamma_2(H)$ and hence $|\gamma_2(H)|\ge 4$. Since $d(H)=3$, $|\Phi(H)|=4$ and thus $Z(H)=\gamma_2(H)=\Phi(H)$. Now $\exp(H/Z(H))=\exp(\gamma_2(H))$ implies that $\gamma_2(H)$ is elementary abelian. As proved in Theorem 4.2, the group $G_{43}$ of fifth family $\Phi_5$ is an extra-special group. Since $\gamma_2(G)$ is family invariant, all the groups in fourth family have elementary abelian center and  all the groups in fifth family have cyclic center. Hence $\mathrm{Aut}_c(G)=\mathrm{Aut}_z(G)$ if and only if  $G\in \Phi_{5}$.
\hfill $\Box$\\

In the next lemma, all the presentations for groups are taken from
[8]. For the sake of brevity, we only write the relators. And for the same reason, a relator of the form $[\alpha_i,\alpha_j]$ is excluded from the presentation if $\alpha_i$ commutes with $\alpha_j$.  For the group $G_{(i,j)}$, $i$ stands for the group number and
$j$ stands for the family number.

\begin{lm}
A non-abelian group $G$ of order $2^6$ is a Camina group of class $2$ if and only if $G\in \Phi_{13}$.
\end{lm}
\noindent {\bf Proof.} Since $\gamma_2(G)=Z(G)$ for a Camina $p$-group of class 2, we need only consider the families $\Phi_9$ to $\Phi_{13}$. These families contain stem groups and the
structure of conjugacy classes for stem groups is an invariant of the isoclinism family. Therefore, it is sufficient to pick one group from each of these families. First consider
$$G=G_{(144,9)}=\langle\alpha_1^{2},\alpha_2^{2},\alpha_3^{2},
\alpha_4^{2},\alpha_5^{2},\alpha_6^{2},[\alpha_4,\alpha_5]
\alpha_1^{-1},[\alpha_4,\alpha_6]\alpha_2^{-1},
[\alpha_5,\alpha_6]\alpha_3^{-1}\rangle.$$
Since $\alpha_1,\alpha_2,\alpha_3\in Z(G)$, $|Z(G)|=|\gamma_2(G)|\geq 8$. Therefore if $x\in G-Z(G)$, then $|x^G|=|[x,G]|\leq 4$ and hence $G$ is not a Camina group. Next consider
$$G=G_{(154,10)}=\langle\alpha_1^{2},\alpha_2^{2},\alpha_4^{2},
\alpha_6^{2},\alpha_3^{2}\alpha_1^{-1},
\alpha_5^{2}\alpha_2^{-1},[\alpha_3,\alpha_4]
\alpha_1^{-1},[\alpha_5,\alpha_6]\alpha_2^{-1}\rangle.$$
Here $\alpha_1,\alpha_2\in Z(G)$, therefore $|Z(G)|=|\gamma_2(G)|\geq 4$.
Since $\alpha_1,\alpha_2,\alpha_3,\alpha_5,\alpha_6\in C_G(\alpha_3)$ and every element of subgroup generated by
$\{\alpha_1,\alpha_2,\alpha_3,\alpha_5,\alpha_6\}$ can be written as $\alpha_1^{i}\alpha_2^{j}\alpha_3^{k}\alpha_5^{l}
\alpha_6^{m},\;0\leq i,j,k,l,m\leq 1$, it follows that
$|C_G(\alpha_3)|=32$. Thus $|\alpha_3^{G}|=|[\alpha_3,G]|=2\neq |\gamma_2(G)|$ and hence $G$ is not a Camina group. Next consider
$$G_{(169,11)}=\langle\alpha_1^{2},\alpha_2^{2},\alpha_3^{2},
\alpha_4^{2},\alpha_6^{2},\alpha_5^{2}\alpha_2^{-1},
[\alpha_4,\alpha_5]
\alpha_1^{-1},[\alpha_3,\alpha_6]\alpha_1^{-1},
[\alpha_5,\alpha_6]\alpha_2^{-1}\rangle.$$
As above, we can show that
$|\alpha_3^{G}|=|[\alpha_3,G]|=2\neq |\gamma_2(G)|\geq 4$ and hence $G_{(169,11)}$ is also not a Camina group. We next consider $G=G_{(180,12)}=$
$$\langle\alpha_1^{2},\alpha_2^{2}\alpha_1^{-1},
\alpha_3^{2}\alpha_1^{-1}, \alpha_4^{2},\alpha_5^{2}\alpha_3^{-1},
\alpha_6^{2}\alpha_4^{-1}, [\alpha_4,\alpha_5]\alpha_1^{-1},
[\alpha_3,\alpha_6]\alpha_1^{-1},
[\alpha_5,\alpha_6]\alpha_2^{-1}\rangle.$$ If $|Z(G)|=16$, then $G/Z(G)$
and hence $\gamma_2(G)$ is elementary abelian, which is not so
because $\alpha_2=[\alpha_5,\alpha_6]\in\gamma_2(G)$ is of order 4.
If $|Z(G)|=8$, then $Z(G)\approx C_2^{2}\times C_2$ or $Z(G)\approx
C_2\times C_2\times C_2$. The first case is not possible by Morigi
Lemma and the second one is ruled out because
$\alpha_2=[\alpha_5,\alpha_6]\in\gamma_2(G)=Z(G)$ is of order 4.
Thus $Z(G)$ is a cyclic group, generated by $\alpha_2$, of order 4.
Therefore $G$ cannot be a Camina group by [14, Lemma 2.1, Theorem 2.2]. Finally consider
$$\begin{array}{ccccc}
G&=&G_{(183,13)}&=&
\langle\alpha_1^{2},\alpha_2^{2},\alpha_3^{2}\alpha_1^{-1},
\alpha_4^{2}\alpha_2^{-1},\alpha_5^{2},\alpha_6^{2},
[\alpha_3,\alpha_5]\alpha_1^{-1},\\
&&&&\\
&&&&[\alpha_4,\alpha_5]\alpha_2^{-1}\alpha_1^{-1},
[\alpha_3,\alpha_6]\alpha_2^{-1}\alpha_1^{-1},
[\alpha_4,\alpha_6]\alpha_2^{-1}\rangle.
\end{array}$$
Observe that $\gamma_2(G)=Z(G)=\langle\alpha_1,\alpha_2\rangle$ and every element $g$ of $G$ is of the form
$z\alpha_3^{m_3}\alpha_4^{m_4}\alpha_5^{m_5}\alpha_6^{m_6}$, where $z\in Z(G)$ and $0\leq m_3,m_4,m_5,m_6\leq 1$.
We claim that $Z(G)\subseteq [g,G]$ for all $g\in G-Z(G)$. First suppose that $m_6=1$. If $m_5=0$, then
$[g,\alpha_4]=[\alpha_5^{m_5}\alpha_6,\alpha_4]=\alpha_2$ where as $[g,\alpha_3]=\alpha_1\alpha_2$. And if $m_5=1$, then $[g,\alpha_4]=\alpha_1$ where as $[g,\alpha_3]=\alpha_2$. Now suppose that
$m_6=0$ and $m_5=1$. Then
$[g,\alpha_3] =\alpha_1$ and $[g,\alpha_4]
=\alpha_1\alpha_2$. We next suppose that $m_5=m_6=0$ and $m_4=1$. If $m_3=0$, then $[g,\alpha_6]=\alpha_2$ where as $[g,\alpha_5]=\alpha_1\alpha_2$. And if $m_3=1$, then $[g,\alpha_6]=\alpha_1$ where as $[g,\alpha_5]=\alpha_2$. Finally suppose that $m_4=m_5=m_6=0$ and $m_3=1$. Then
$[g,\alpha_5]=\alpha_1$ and $[g,\alpha_6]=\alpha_1\alpha_2$. In every possibility, $Z(G)\subseteq [g,G]$ for all $g\in G-Z(G)$. This proves the claim and the lemma.\hfill $\Box$

\begin{thm}
Let $G$ be a non-abelian group of order $2^6$. Then $\mathrm{Aut}_c(G)=\mathrm{Aut}_z(G)$ if and only if $G\in \Phi_{12}$ or $\Phi_{13}$.
\end{thm}
\noindent {\bf Proof.} First suppose that
$\mathrm{Aut}_c(G)=\mathrm{Aut}_z(G)$. By Theorem 4.4, either $G$ is
a Camina group of class 2 or $\gamma_2(G)=Z(G)$ and $Z(G)$ is
cyclic. If $G$ is a Camina group of class 2, then $G\in
\Phi_{13}$ by above lemma. Assume $\gamma_2(G)=Z(G)$ and $Z(G)$ is
cyclic. There are only five isoclinism families $viz.$ $\Phi_9$ to
$\Phi_{13}$ which consist of groups $G$ such that
$\gamma_2(G)=Z(G)$. Consider $G_{144}\in \Phi_9$. By above lemma
$|Z(G_{144})|\ge 8$ and by [16], $d(G_{144})=3$. Thus
$Z(G_{144})=\gamma_2(G_{144})=\Phi(G_{144})$. Now
$\exp(G_{144}/Z(G_{144}))=\exp(\gamma_2(G_{144}))$ implies that
$Z(G_{144})$ is elementary abelian of order 8. Next consider
$G_{154}\in \Phi_{10}$. By above lemma $|Z(G_{154})|\ge 4$ and by
[16], $d(G_{154})=4$. Thus
$Z(G_{154})=\gamma_2(G_{154})=\Phi(G_{154})$. Then
$\exp(G_{154}/Z(G_{154}))=\exp(\gamma_2(G_{154}))$ implies that
$Z(G_{154})$ is elementary abelian of order 4. Next consider
$G_{169}\in \Phi_{11}$. By above lemma $|Z(G_{169})|\ge 4$ and by
[16], $d(G_{169})=4$. Thus
$Z(G_{169})=\gamma_2(G_{169})=\Phi(G_{169})$. Therefore
$\exp(G_{169}/Z(G_{169}))=\exp(\gamma_2(G_{169}))$ implies that
$Z(G_{169})$ is elementary abelian of order 4. Finally, if $G\in
\Phi_{12}$, then $Z(G)$ is cyclic by Lemma 5.5.\hfill $\Box$

\vspace{.2in} \noindent {\bf Acknowledgements.} Authors are thankful to the referee for his useful comments and suggestions. The research of the
first author is supported by Council of Scientific and Industrial
Research, Government of India. The same is gratefully acknowledged.

\end{document}